\documentclass[12pt]{article}
\usepackage{amsfonts}

\usepackage{mathrsfs}
\usepackage{amscd}
\usepackage{amsmath,amsfonts,amssymb,amscd}
\usepackage{indentfirst,graphics,epsfig,psfrag}
\input{epsf}

\usepackage{amsmath,amssymb,amsthm, amscd, amsfonts, graphicx,ifpdf}
\usepackage{lineno}

\newtheorem{thm}{Theorem}[section]

\newtheorem{lem}{Lemma}[section]

\def\qed{\nopagebreak\hfill{\rule{4pt}{7pt}}
\medbreak}

\def\pf{\noindent {\it Proof.} }

\title{\bf The generalized connectivity of\\ complete bipartite graphs\footnote{Supported by NSFC.}}

\author{
\small Shasha Li, Wei Li, Xueliang Li\\
\small Center for Combinatorics and LPMC-TJKLC\\
\small Nankai University, Tianjin 300071, China.\\
\small  Email: lss@cfc.nankai.edu.cn, liwei@cfc.nankai.edu.cn, lxl@nankai.edu.cn\\
}
\date{}
\begin{document}

\maketitle

\begin{abstract}
Let $G$ be a nontrivial connected graph of order $n$, and $k$ an
integer with $2\leq k\leq n$. For a set $S$ of $k$ vertices of $G$,
let $\kappa (S)$ denote the maximum number $\ell$ of edge-disjoint
trees $T_1,T_2,\ldots,T_\ell$ in $G$ such that $V(T_i)\cap V(T_j)=S$
for every pair $i,j$ of distinct integers with $1\leq i,j\leq \ell$.
Chartrand et al. generalized the concept of connectivity as follows:
The $k$-$connectivity$, denoted by $\kappa_k(G)$, of $G$ is defined
by $\kappa_k(G)=$min$\{\kappa(S)\}$, where the minimum is taken over
all $k$-subsets $S$ of $V(G)$. Thus $\kappa_2(G)=\kappa(G)$, where
$\kappa(G)$ is the connectivity of $G$. Moreover, $\kappa_{n}(G)$ is
the maximum number of edge-disjoint spanning trees of $G$.

This paper mainly focus on the $k$-connectivity of complete
bipartite graphs $K_{a,b}$. First, we obtain the number of
edge-disjoint spanning trees of $K_{a,b}$, which is
$\lfloor\frac{ab}{a+b-1}\rfloor$, and specifically give the
$\lfloor\frac{ab}{a+b-1}\rfloor$ edge-disjoint spanning trees. Then
based on this result, we get the $k$-connectivity of $K_{a,b}$ for
all $2\leq k \leq a+b$. Namely, if $k>b-a+2$ and $a-b+k$ is odd then
$\kappa_{k}(K_{a,b})=\frac{a+b-k+1}{2}+\lfloor\frac{(a-b+k-1)(b-a+k-1)}{4(k-1)}\rfloor,$
if $k>b-a+2$ and $a-b+k$ is even then
$\kappa_{k}(K_{a,b})=\frac{a+b-k}{2}+\lfloor\frac{(a-b+k)(b-a+k)}{4(k-1)}\rfloor,$
and if $k\leq b-a+2$ then $\kappa_{k}(K_{a,b})=a. $
\\[3mm]
\end{abstract}

\section{Introduction}

We follow the terminology and notation of \cite{Bondy}. As usual,
denote by $K_{a,b}$ the complete bipartite graph with bipartition of
sizes $a$ and $b$. The $connectivity$ $\kappa(G)$ of a graph $G$ is
defined as the minimum cardinality of a set $Q$ of vertices of $G$
such that $G-Q$ is disconnected or trivial. A well-known theorem of
Whitney \cite{Whitney} provides an equivalent definition of the
connectivity. For each $2$-subset $S=\{u,v\}$ of vertices of $G$,
let $\kappa(S)$ denote the maximum number of internally disjoint
$uv$-paths in $G$. Then $\kappa(G)=$min$\{\kappa(S)\}$, where the
minimum is taken over all $2$-subsets $S$ of $V(G)$.

In \cite{Chartrand}, the authors generalized the concept of
connectivity. Let $G$ be a nontrivial connected graph of order $n$,
and $k$ an integer with $2\leq k\leq n$. For a set $S$ of $k$
vertices of $G$, let $\kappa (S)$ denote the maximum number $\ell$
of edge-disjoint trees $T_1,T_2,\ldots,T_\ell$ in $G$ such that
$V(T_i)\cap V(T_j)=S$ for every pair $i,j$ of distinct integers with
$1\leq i,j\leq \ell$ (Note that the trees are vertex-disjoint in
$G\backslash S$). A collection $\{T_1,T_2,\ldots,T_\ell \}$ of trees
in $G$ with this property is called an {\it internally disjoint set
of trees connecting $S$}. The $k$-$connectivity$, denoted by
$\kappa_k(G)$, of $G$ is then defined as
$\kappa_k(G)=$min$\{\kappa(S)\}$, where the minimum is taken over
all $k$-subsets $S$ of $V(G)$. Thus, $\kappa_2(G)=\kappa(G)$ and
$\kappa_{n}(G)$ is the maximum number of edge-disjoint spanning
trees of $G$.

In \cite{LLZ}, the authors focused on the investigation of
$\kappa_{3}(G)$ and mainly studied the relationship between the
$2$-connectivity and the $3$-connectivity of a graph. They gave
sharp upper and lower bounds for $\kappa_{3}(G)$ for general graphs
$G$, and showed that if $G$ is a connected planar graph, then
$\kappa(G)-1\leq \kappa_{3}(G)\leq\kappa(G)$. Moreover, they studied
the algorithmic aspects for $\kappa_{3}(G)$ and gave an algorithm to
determine $\kappa_{3}(G)$ for a general graph $G$.

Chartrand et al. in \cite{Chartrand} proved that if $G$ is the
complete $3$-partite graph $K_{3,4,5}$, then $\kappa_3(G)=6$. They
also gave a general result for the complete graph $K_n$:
\begin{thm}\label{thm0}
For every two integers $n$ and $k$ with $2\leq k\leq n$,
\begin{center}
$\kappa_k(K_n)=n-\lceil k/2\rceil$.
\end{center}
\end{thm}

In this paper, we turn to complete bipartite graphs $K_{a,b}$.
First, we give the number of edge-disjoint spanning trees of
$K_{a,b}$, namely $\kappa_{a+b}(K_{a,b})$.
\begin{thm}\label{thm1}
For every two integers $a$ and $b$,
$$
\kappa_{a+b}(K_{a,b}) = \lfloor\frac{ab}{a+b-1}\rfloor.
$$
\end{thm}
Actually, we specifically give the $\lfloor\frac{ab}{a+b-1}\rfloor$
edge-disjoint spanning trees of $K_{a,b}$. Then based on Theorem
\ref{thm1}, we obtain the $k$-connectivity of $K_{a,b}$ for all
$2\leq k \leq a+b$.

\section{Proof of Theorem \ref{thm1}}

Since $K_{a,b}$ contains $ab$ edges and a spanning tree needs
$a+b-1$ edges, the number of edge-disjoint spanning trees of
$K_{a,b}$ is at most $\lfloor\frac{ab}{a+b-1}\rfloor$, namely,
$\kappa_{a+b}(K_{a,b}) \leq \lfloor\frac{ab}{a+b-1}\rfloor$. Thus,
it suffices to prove that $\kappa_{a+b}(K_{a,b}) \geq
\lfloor\frac{ab}{a+b-1}\rfloor$. To this end, we want to find all
the $\lfloor\frac{ab}{a+b-1}\rfloor$ edge-disjoint spanning trees.

Let $X=\{x_{1}, x_{2}, \ldots , x_{a}\}$ and $Y=\{y_{1}, y_{2},
\ldots , y_{b}\}$ be the bipartition of $K_{a,b}$. Without loss of
generality, we may assume that $a\leq b$.

We will express the spanning trees by adjacency-degree lists. To be
specific, the fist spanning tree $T_{1}$ we find can be represented
by an adjacency-degree list as follows:

\begin{tabular}{l|l|l}
vertex & neighbors & degree \\
\hline
$x_{1}$ & $y_{1},\ y_{2},\ \ldots,\ y_{d_{1}}$ & $d_{1}$\\
$x_{2}$ & $y_{d_{1}},\ y_{d_{1}+1},\ \ldots,\ y_{d_{1}+d_{2}-1}$ & $d_{2}$\\
$x_{3}$ & $y_{d_{1}+d_{2}-1},\ y_{d_{1}+d_{2}},\ \ldots,\ y_{d_{1}+d_{2}+d_{3}-2}$ & $d_{3}$\\
$\ldots$ & $\ldots$ & $\ldots$\\
$x_{j}$ & $y_{d_{1}+d_{2}+\cdots+d_{j-1}-(j-2)},\ y_{d_{1}+d_{2}+\cdots+d_{j-1}-(j-2)+1},\ \ldots,\
y_{d_{1}+d_{2}+\cdots+d_{j}-(j-1)}$ & $d_{j}$\\
$\ldots$ & $\ldots$ & $\ldots$\\
$x_{a}$ & $y_{d_{1}+d_{2}+\cdots+d_{a-1}-(a-2)},\ y_{d_{1}+d_{2}+\cdots+d_{a-1}-(a-2)+1},\ \ldots,\
y_{d_{1}+d_{2}+\cdots+d_{a}-(a-1)}$ & $d_{a}$
\end{tabular}

\noindent where $d_{j}$ denotes the degree of $x_{j}$ in $T_{1}$,
and $d_{1} + d_{2} + \cdots + d_{a}= a+b-1$.

To simplify the subscript, we denote $i_{0}=1$, $i_{1}=d_{1}$, $i_{2}=d_{1}+d_{2}-1$, $\ldots$,
$i_{j}=d_{1}+d_{2}+\cdots+d_{j}-(j-1)$, $\ldots$, $i_{a}=d_{1}+d_{2}+\cdots+d_{a}-(a-1)=b$. Note that,
$i_{j}-i_{j-1}=d_{j}-1$. So the adjacency-degree list of $T_{1}$ can be simplified as follows:

$T_{1}$\begin{tabular}{l|l|l}
vertex & neighbors & degree\\
\hline
$x_{1}$ & $y_{i_{0}},\ y_{i_{0}+1},\ \ldots,\ y_{i_{1}}$ & $d_{1}$\\
$x_{2}$ & $y_{i_{1}},\ y_{i_{1}+1},\ \ldots,\ y_{i_{2}}$ & $d_{2}$\\
$x_{3}$ & $y_{i_{2}},\ y_{i_{2}+1},\ \ldots,\ y_{i_{3}}$ & $d_{3}$\\
$\ldots$ & $\ldots$ & $\ldots$\\
$x_{j}$ & $y_{i_{j-1}},\ y_{i_{j-1}+1},\ \ldots,\ y_{i_{j}}$ & $d_{j}$\\
$\ldots$ & $\ldots$ & $\ldots$\\
$x_{a}$ & $y_{i_{a-1}},\ y_{i_{a-1}+1},\ \ldots,\ y_{i_{a}}$ & $d_{a}$
\end{tabular}

Then we can list the second spanning trees we find. Here and in what follows, for a vertex $y_{j}$, if $j>b$, $y_{j}$
denotes $y_{j-b}$, for a subscript $i_{j}$, if $j>a$, $y_{i_{j}}$ denotes $y_{i_{j-a}}$, and for degree $d_{j}$, if
$j>a$, $d_{j}$ denotes $d_{j-a}$.

$T_{2}$
\begin{tabular}{l|l|l}
vertex & neighbors & degree\\
\hline
$x_{1}$ & $y_{i_{1}+1},\ y_{i_{1}+2},\ \ldots,\ y_{i_{2}+1}$ & $d_{2}$\\
$x_{2}$ & $y_{i_{2}+1},\ y_{i_{2}+2},\ \ldots,\ y_{i_{3}+1}$ & $d_{3}$\\
$x_{3}$ & $y_{i_{3}+1},\ y_{i_{3}+2},\ \ldots,\ y_{i_{4}+1}$ & $d_{4}$\\
$\ldots$ & $\ldots$ & $\ldots$\\
$x_{j}$ & $y_{i_{j}+1},\ y_{i_{j}+2},\ \ldots,\ y_{i_{j+1}+1}$ & $d_{j+1}$\\
$\ldots$ & $\ldots$ & $\ldots$\\
$x_{a}$ & $y_{i_{a}+1},\ y_{i_{a}+2},\ \ldots,\ y_{i_{a+1}}$ & $d_{1}$
\end{tabular}

From the lists, we can see that $T_{2}$ and $T_{1}$ are edge-disjoint, if and only if for every vertex $x_{j}$,
$d_{j}+d_{j+1}\leq b$. If $T_{2}$ and $T_{1}$ are edge-disjoint, then we continue to list $T_{3}$.

$T_{3}$
\begin{tabular}{l|l|l}
vertex & neighbors & degree\\
\hline
$x_{1}$ & $y_{i_{2}+2},\ y_{i_{2}+3},\ \ldots,\ y_{i_{3}+2}$ & $d_{3}$\\
$x_{2}$ & $y_{i_{3}+2},\ y_{i_{3}+3},\ \ldots,\ y_{i_{4}+2}$ & $d_{4}$\\
$x_{3}$ & $y_{i_{4}+2},\ y_{i_{4}+3},\ \ldots,\ y_{i_{5}+2}$ & $d_{5}$\\
$\ldots$ & $\ldots$ & $\ldots$\\
$x_{j}$ & $y_{i_{j+1}+2},\ y_{i_{j+1}+3},\ \ldots,\ y_{i_{j+2}+2}$ & $d_{j+2}$\\
$\ldots$ & $\ldots$ & $\ldots$\\
$x_{a}$ & $y_{i_{a+1}+2},\ y_{i_{a+1}+3},\ \ldots,\ y_{i_{a+2}+1}$ & $d_{2}$
\end{tabular}

From the lists, we can see that $T_{3}$ and $T_{1}$, $T_{2}$ are edge-disjoint, if and only if for every vertex
$x_{j}$, $d_{j}+d_{j+1}+d_{j+2}\leq b$. If $T_{3}$ and $T_{1}$, $T_{2}$ are edge-disjoint, then we continue to list
$T_{4}$. Continuing the procedure, our goal is to find the maximum $t$, such that $T_{t}$ and $T_{1}, T_{2}, \ldots,
T_{t-1}$ are edge-disjoint.

$T_{t}$
\begin{tabular}{l|l|l}
vertex & neighbors & degree\\
\hline
$x_{1}$ & $y_{i_{t-1}+(t-1)},\ y_{i_{t-1}+t},\ \ldots,\ y_{i_{t}+(t-1)}$ & $d_{t}$\\
$x_{2}$ & $y_{i_{t}+(t-1)},\ y_{i_{t}+t},\ \ldots,\ y_{i_{t+1}+(t-1)}$ & $d_{t+1}$\\
$x_{3}$ & $y_{i_{t+1}+(t-1)},\ y_{i_{t+1}+t},\ \ldots,\ y_{i_{t+2}+(t-1)}$ & $d_{t+2}$\\
$\ldots$ & $\ldots$\\
$x_{j}$ & $y_{i_{j+t-2}+(t-1)},\ y_{i_{j+t-2}+t},\ \ldots,\ y_{i_{j+t-1}+(t-1)}$ & $d_{t+j-1}$\\
$\ldots$ & $\ldots$\\
$x_{a}$ & $y_{i_{a+t-2}+(t-1)},\ y_{i_{a+t-2}+t},\ \ldots,\ y_{i_{a+t-1}+(t-2)}$ & $d_{t-1}$
\end{tabular}

\noindent That is, we want to find the maximum $t$, such that
$d_{j}+d_{j+1}+\cdots+d_{j+t-1}\leq b$, for any $1\leq j\leq a$.

Let $D_{j}^{t}=d_{j} + d_{j+1} + \cdots + d_{j+t-1}$. It can be
observed that $D_{j}^{t}=D_{j+1}^{t}$ if and only if
$d_{j}=d_{j+t}$. Consider the numbers $1, t+1, 2t+1, \ldots,
(a-1)t+1$, where addition is carried out by modula $a$.

\textbf{Case 1.} $1, t+1, 2t+1, \ldots, (a-1)t+1$ are pairwise distinct.

Then we can assign the values to $d_{j}$ as follows:

Let $a+b-1=ka+c$, where $k, c$ are integers, and $0\leq c\leq a-1$.
Then $a+b-1=(k+1)c+k(a-c)$. If $c=0$, let $d_{j}=k$, for all $1\leq
j \leq a$. If $c\neq 0$, let $d_{it+1}=k+1$, for all $0\leq i\leq
c-1$, and let other $d_{j}=k$.

\textbf{Case 2.} Some of the numbers $1, t+1, 2t+1, \ldots ,
(a-1)t+1$ are equal.

Without loss of generality, suppose $jt+1$ is the first number that
equals a number $it+1$ before it, namely, $jt+1=it+1 \ (mod\ a)$,
where $j >i$. Then $(j-i)t+1=1 \ (mod\ a)$. Since $jt+1$ is the
first number that equals a number before it, we can get $i=0$. Thus,
$1, t+1, 2t+1, \ldots, (j-1)t+1$ are pairwise distinct.

\textbf{Claim 1.} $it+1 \neq 2 \ (mod\ a)$, for any integer $i$.

If $it+1=2 \ (mod\ a)$, then we have $it=1 \ (mod\ a)$. Thus we have
$$
\begin{array}{cccccc}
it&+&1&=&2&(mod\ a)\\
2it&+&1&=&3&(mod\ a)\\
\multicolumn{6}{c}{\dotfill} \\
(a-1)it&+&1&=&a&(mod\ a)
\end{array}
$$
So there are $a$ distinct numbers in $\{1, it+1, 2it+1, \ldots ,
(a-1)it+1\}$. On the other hand, since $jt+1=1 \ (mod\ a)$, there
are at most $j\leq a-1$ distinct numbers in $\{ut+1, \text{$u$ is an
integer}\}\supset\{1, it+1, 2it+1, \ldots , (a-1)it+1\}$, a
contradiction. Thus,  $it+1 \neq 2 \ (mod\ a)$ for any integer $i$.

\textbf{Claim 2. } $2, t+2, 2t+2, \ldots , (j-1)t+2$ are pairwise distinct.

If $j_{1}t+2=j_{2}t+2 \ (mod\ a)$, where $0\leq j_{1} < j_{2} \leq
j-1$, then $j_{1}t+1=j_{2}t+1 \ (mod\ a)$. But $1, t+1, 2t+1,
\ldots, (j-1)t+1$ are pairwise distinct, a contradiction. Thus, $2,
t+2, 2t+2, \ldots , (j-1)t+2$ are pairwise distinct.

\textbf{Claim 3.} $\{1, t+1, 2t+1, \ldots, (j-1)t+1\}\cap\{2, t+2, 2t+2, \ldots , (j-1)t+2\}=\emptyset$.

If $i_{1}t+1=i_{2}t+2 \ (mod\ a)$, then $(i_{1}-i_{2})t+1=2 \ (mod\
a)$, but $it+1 \neq 2 \ (mod\ a)$ for any integer $i$, a
contradiction by Claim 1. Thus, $1, t+1, 2t+1, \ldots, (j-1)t+1, 2,
t+2, 2t+2, \ldots , (j-1)t+2$ are pairwise distinct.

Now, if $2=\frac{a}{j}$, then we have already ordered all numbers of
$\{1,\ \ldots,\ a\}$. Else if $2<\frac{a}{j}$, we will prove that
$1+it\neq3 \ (mod\ a)$ and $2+it\neq3 \ (mod\ a)$ for any integer
$i$.

\textbf{Claim 4.} If $2<\frac{a}{j}$, then $1+it\neq3 \ (mod\ a)$
and $2+it\neq3 \ (mod\ a)$ for any integer $i$.

If $2+it=3 \ (mod\ a)$, then $1+it=2 \ (mod\ a)$, a contradiction by
Claim 1. If $1+it=3 \ (mod\ a)$, then we have $it=2 \ (mod\ a)$.
Thus we have
$$
\begin{array}{cccccc}
it&+&1&=&3&(mod\ a)\\
it&+&2&=&4&(mod\ a)\\
2it&+&1&=&5&(mod\ a)\\
2it&+&2&=&6&(mod\ a)\\
\multicolumn{6}{c}{\dotfill} \\
\frac{a-2}{2}it&+&1&=&a-1& (mod\ a) \ (\text{for $a$ even})\\
\frac{a-3}{2}it&+&2&=&a-1& (mod\ a) \ (\text{for $a$ odd})\\
\frac{a-2}{2}it&+&2&=&a& (mod\ a) \ (\text{for $a$ even})\\
\frac{a-1}{2}it&+&1&=&a& (mod\ a) \ (\text{for $a$ odd})
\end{array}
$$
So there are at least $a$ distinct numbers in $\{1, it+1, 2it+1,
\ldots , \lceil\frac{a}{2}\rceil it+1, 2, it+2, 2it+2, \ldots ,
\lceil\frac{a}{2}\rceil it+2\}$. On the other hand, since $jt+1=1 \
(mod\ a)$ and $j \leq a-1$, there are at most $2j<a$ distinct
numbers in $\{ut+1, \text{$u$ is an integer}\}\cup\{vt+2, \text{$v$
is an integer}\}\supset\{1, it+1, 2it+1, \ldots ,
\lceil\frac{a}{2}\rceil it+1, 2, it+2, 2it+2, \ldots ,
\lceil\frac{a}{2}\rceil it+2\}$, a contradiction. Hence, if
$2<\frac{a}{j}$, then $1+it\neq3 \ (mod\ a)$ and $2+it\neq3 \ (mod\
a)$ for any integer $i$.

Similarly, we can prove that $r+it\neq s \ (mod\ a)$ for $1\leq
r<s\leq\frac{a}{j}$. Thus we can get the following claim:

\textbf{Claim 5.} $1, t+1, 2t+1, \ldots, (j-1)t+1, 2, t+2, 2t+2, \ldots , (j-1)t+2, \ldots, \frac{a}{j},
t+\frac{a}{j}, 2t+\frac{a}{j}, \ldots , (j-1)t+\frac{a}{j}$ are pairwise distinct. And hence $\{1, t+1, 2t+1, \ldots,
(j-1)t+1\}\cup\{2, t+2, 2t+2, \ldots , (j-1)t+2\}\cup\cdots\cup\{\frac{a}{j}, t+\frac{a}{j}, 2t+\frac{a}{j}, \ldots ,
(j-1)t+\frac{a}{j}\}=\{1, 2, \ldots, a\}$.

The proof is similar to those of Claims 2, 3 and 4. We thus have
ordered $\{1, 2, \ldots, a\}$ by $1, t+1, 2t+1, \ldots, (j-1)t+1, 2,
t+2, 2t+2, \ldots , (j-1)t+2, \ldots, \frac{a}{j}, t+\frac{a}{j},
2t+\frac{a}{j}, \ldots , (j-1)t+\frac{a}{j}$. Let $a+b-1=ka+c$,
where $k, c$ are integers, and $0\leq c\leq a-1$. Then
$a+b-1=(k+1)c+k(a-c)$.

Now, we can assign the values of $d_{j}$ as follows: If $c=0$, let
$d_{j}=k$ for all $1\leq j \leq a$. If $c\neq 0$, for the first $c$
numbers of our ordering, if $d_{j}$ uses one of them as subscript,
then $d_{j}=k+1$, else $d_{j}=k$.

Next, we will show that, in either case, $\mid
D_{i}^{t}-D_{j}^{t}\mid \leq 1$ for any integers $1\leq i, j \leq a$
and $t>0$.

If $c=0$, $d_{j}=k$ for all $1\leq j \leq a$, then
$D_{i}^{t}=D_{j}^{t}$ for any integers $1\leq i, j \leq a$. The
assertion is certainly true. So we may assume that $c\neq 0$. For
Case 1, we construct a weighted cycle: $C=v_{1}v_{2}\ldots
v_{a}v_{1}$ and $w(v_{i})= d_{(i-1)t+1}$, where $v_{i}$ corresponds
to vertex $x_{(i-1)t+1}$, $1\leq i \leq a$.

According to the assignment,
$$
w(v_{1})=w(v_{2})=\cdots =w(v_{c})=k+1,
$$
and
$$
w(v_{c+1})=w(v_{c+2})=\cdots =w(v_{a})=k.
$$
Since $D_{i}^{t}=D_{i+1}^{t}$ if and only if $d_{i}=d_{i+t}$, then
$D_{(i-1)t+1}^{t}=D_{(i-1)t+1+1}^{t}$ if and only if
$w(v_{i})=w(v_{i+1}).$ Similarly,
$D_{(i-1)t+1}^{t}=D_{(i-1)t+1+1}^{t}+1$ if and only if
$w(v_{i})=w(v_{i+1})+1,$ and $
D_{(i-1)t+1}^{t}=D_{(i-1)t+1+1}^{t}-1$ if and only if
$w(v_{i})=w(v_{i+1})-1. $ We know that $w(v_{c})=w(v_{c+1})+1$ and
$w(v_{a})=w(v_{1})-1$. For simplicity, let $(c-1)t+1=\alpha \ (mod\
a)$, $(a-1)t+1=\beta \ (mod\ a)$, that is, $v_{c}$ corresponds to
$x_{\alpha}$ and $v_{a}$ corresponds to $x_{\beta}$, and by the
hypothesis, $\alpha\neq\beta$.

If $\alpha < \beta$, then
$$
D_{1}^{t}=D_{2}^{t}=\cdots =D_{\alpha}^{t}=D_{\alpha+1}^{t}+1=D_{\alpha+2}^{t}+1=\cdots
=D_{\beta}^{t}+1=D_{\beta+1}^{t}=D_{\beta+2}^{t}=\cdots=D_{a}^{t}.
$$
If $\alpha > \beta$, then
$$
D_{1}^{t}=D_{2}^{t}=\cdots =D_{\beta}^{t}=D_{\beta+1}^{t}-1=D_{\beta+2}^{t}-1=\cdots
=D_{\alpha}^{t}-1=D_{\alpha+1}^{t}=D_{\alpha+2}^{t}=\cdots=D_{a}^{t}.
$$
In any case, we have $\mid D_{i}^{t}-D_{j}^{t}\mid \leq 1$ for any
integers $1\leq i, j \leq a$ and $t>0$.

For Case 2, we construct $\frac{a}{j}$ weighted cycles.
$C_{i}=v_{i_{1}}v_{i_{2}}\ldots v_{i_{j}}v_{i_{1}}$, $1 \leq i \leq
\frac{a}{j}$, and $w(v_{i_{r}})= d_{(r-1)t+i}$, where $v_{i_{r}}$
corresponds to vertex $x_{(r-1)t+i}$, $1\leq r \leq j$. By the
assignment, there is at most one cycle in which the vertices have
two distinct weights. If such cycle does not exist, clearly, we have
$D_{1}^{t}=D_{2}^{t}=\cdots =D_{a}^{t}$. So we may assume that for
some cycle $C_{s}$, $w(v_{s_{\gamma}})=w(v_{s_{\gamma+1}})+1$ and
$w(v_{s_{j}})=w(v_{s_{1}})-1$. Similar to the proof of Case 1, we
can get that $\mid D_{i}^{t}-D_{j}^{t}\mid \leq 1$ for any integers
$1\leq i, j \leq a$ and $t>0$.

Then, we can show that, with the assignment we can get $t\geq \lfloor\frac{ab}{a+b-1}\rfloor$.

Let $t'=\lfloor\frac{ab}{a+b-1}\rfloor$. And let
$$
\begin{array}{ccccccccc}
    D_{1}^{t'} &=& d_{1} &+& d_{2} &+& \cdots &+& d_{t'}\\
    D_{2}^{t'} &=& d_{2} &+& d_{3} &+& \cdots &+& d_{t'+1}\\
    \multicolumn{9}{c}{\dotfill} \\
    D_{j}^{t'} &=& d_{j} &+& d_{j+1} &+& \cdots &+& d_{j+t'-1}\\
    \multicolumn{9}{c}{\dotfill} \\
    D_{a}^{t'} &=& d_{a} &+& d_{1} &+& \cdots &+& d_{t'-1}
\end{array}
$$
we have
$D_{1}^{t'}+D_{2}^{t'}+\cdots + D_{a}^{t'} =t'(d_{1}+d_{2}+\cdots+d_{a})= t'(a+b-1)$

It follows from $\mid D_{i}^{t}-D_{j}^{t}\mid \leq 1$, for any integers $1\leq i, j \leq a$ and $t>0$, that
$$
D_{j}^{t'}\leq \lceil\frac{t'(a+b-1)}{a}\rceil < \frac{t'(a+b-1)}{a} + 1 \leq \frac{ab}{a+b-1}\frac{a+b-1}{a}+1= b+1
$$
The third inequality holds since $t' =
\lfloor\frac{ab}{a+b-1}\rfloor \leq \frac{ab}{a+b-1}$. Since
$D_{j}^{t'}$ is an integer, we have $D_{j}^{t'}\leq b$ for all
$1\leq j \leq a$. Since $t$ is the maximum integer such that
$D_{j}^{t}=d_{j}+d_{j+1}+\cdots+d_{j+t-1}\leq b$ for any $1\leq j
\leq a$, then $t\geq t'=\lfloor\frac{ab}{a+b-1}\rfloor$. So we can
find at least $\lfloor\frac{ab}{a+b-1}\rfloor$ edge-disjoint
spanning trees of $K_{a,b}$. And hence
$\kappa_{a+b}(K_{a,b})\geq\lfloor\frac{ab}{a+b-1}\rfloor$. So we
have proved that
$\kappa_{a+b}(K_{a,b})=\lfloor\frac{ab}{a+b-1}\rfloor$. \qed

\section{The $k$-connectivity of complete bipartite graphs}

Next, we will calculate $\kappa_{k}(K_{a,b})$, for $2\leq k \leq a+b$.

Recall that $\kappa_{k}(G) = \min\{\kappa(S)\}$, where the minimum
is taken over all $k$-element subsets $S$ of $V(G)$. Denote by
$K_{a,\ b}$ a complete bipartite graph with bipartition $X=\{x_{1},
x_{2}, \ldots, x_{a}\}$ and $Y=\{y_{1}, y_{2}, \ldots, y_{b}\}$,
where $a\leq b$. Actually, all vertices in $X$ are equivalent and
all vertices in $Y$ are equivalent. So instead of considering all
$k$-element subsets $S$ of $V(G)$, we can restrict our attention to
the subsets $S_{i}$, for $0 \leq i \leq k$, where $S_{i}$ is an
$k$-element subsets of $V(G)$ such that $S_{i}\cap X = \{x_{1},
x_{2}, \ldots, x_{i}\}$, $S_{i}\cap Y = \{y_{1}, y_{2}, \ldots,
y_{k-i}\}$, $1 \leq i \leq k$ and $S_{0}\cap X =\emptyset$,
$S_{0}\cap Y =\{y_{1}, y_{2}, \ldots, y_{k}\}$. Notice that, if
$i>a$ or $k-i > b$ then $S_{i}$ does not exist, and if $k > b$ then
$S_{0}$ does not exist. So, we need only to consider $S_{i}$ for
$\max\{0,\ k-b\} \leq i \leq \min\{a,\ k\}$.

Now, let $A$ be a maximum set of internally disjoint trees
connecting $S_{i}$. Let $\mathfrak{A}_{0}$ be the set of trees
connecting $S_{i}$ whose vertex set is $S_{i}$, let
$\mathfrak{A}_{1}$ be the set of trees connecting $S_{i}$ whose
vertex set is $S_{i}\cup \{u\}$, where $u \notin S_{i}$ and let
$\mathfrak{A}_{2}$ be the set of trees connecting $S_{i}$ whose
vertex set is $S_{i}\cup \{u, v\}$, where $u, v \notin S_{i}$ and
they belong to distinct partitions.

\begin{lem}\label{lem1}
Let $A$ be a maximum set of internally disjoint trees connecting
$S_{i}$. Then we can always find a set $A'$ of internally disjoint
trees connecting $S_{i}$, such that $\mid A \mid = \mid A' \mid$ and
$A' \subset
\mathfrak{A}_{0}\cup\mathfrak{A}_{1}\cup\mathfrak{A}_{2}$.
\end{lem}
\pf If there is a tree $T^{0}$ in $A$ whose vertex set
$V(T^{0})\supseteq\{u_{1}, u_{2}\}$, where $u_{1}, u_{2}\notin
S_{i}$ and $u_{1}, u_{2}$ belong to the same partition, then we can
connect all neighbors of $u_{2}$ to $u_{1}$ by some new edges and
delete $u_{2}$ and the multiple edges (if exist). Obviously, the new
graph we obtain is still a tree $T'$ that connect $S_{i}$. Since
$V(T_{m})\cap V(T_{n}) = S_{i}$ for every pair of trees in $A$,
other trees in $A$ will not contain $u_{1}$, including the edges
incident with $u_{1}$. So for all trees $T_{n}$ in $A$ other than
$T^{0}$, $V(T')\cap V(T_{n}) = S_{i}$ and $E(T')\cap E(T_{n})
=\emptyset$. Moreover, $T'$ has less vertices which are not in
$S_{i}$ than $T^{0}$. Repeat this process, until we get a tree $T
\in \mathfrak{A}_{0}\cup \mathfrak{A}_{1}\cup \mathfrak{A}_{2}$.
Replace $A$ by $A^{1}=A\setminus\{T^{0}\}\cup\{T\}$, and then
$A^{1}$ contains less trees that are not in $\mathfrak{A}_{0}\cup
\mathfrak{A}_{1}\cup \mathfrak{A}_{2}$ than $A$. Repeating the
process, we can get a series of sets $A^{0}, A^{1}, \ldots, A^{t}$,
such that $A^{0}=A$, $A^{t}=A'$, and $A^{j}$ contains less trees not
in $\mathfrak{A}_{0}\cup \mathfrak{A}_{1}\cup \mathfrak{A}_{2}$ than
$A^{j-1}$ for $1\leq j \leq t$, where all $A^{s}$ are sets of
internally disjoint trees connecting $S_{i}$ for $0\leq s \leq t$,
and $\mid A^{0}\mid=\cdots=\mid A^{t}\mid$. So we finally get the
set $A'\subset \mathfrak{A}_{0}\cup \mathfrak{A}_{1}\cup
\mathfrak{A}_{2}$ which has the same cardinality as $A$.\qed

So, we can assume that the maximum set $A$ of internally disjoint
trees connecting $S_{i}$ is contained in
$\mathfrak{A}_{0}\cup\mathfrak{A}_{1}\cup\mathfrak{A}_{2}$.

Next, we will define the standard structure of trees in
$\mathfrak{A}_{0}$, $\mathfrak{A}_{1}$ and $\mathfrak{A}_{2}$,
respectively.

Every tree in $\mathfrak{A}_{0}$ is of standard structure. A tree
$T$ in $\mathfrak{A}_{1}$ with vertex set $V(T)=S_{i}\cup\{u\}$,
where $u\in X\setminus S_{i}$, is of standard structure, if $u$ is
adjacent to every vertex in $S_{i}\cap Y$, and every vertex in
$S_{i}\cap X$ has degree $1$. A tree $T$ in $\mathfrak{A}_{1}$ with
vertex set $V(T)=S_{i}\cup\{v\}$, where $v\in Y\setminus S_{i}$, is
of standard structure, if $v$ is adjacent to every vertex in
$S_{i}\cap X$, and every vertex in $S_{i}\cap Y$ has degree $1$. A
tree $T$ in $\mathfrak{A}_{2}$ with vertex set $V(T)=S_{i}\cup\{u,
v\}$, where $u\in X\setminus S_{i}$ and $v\in Y\setminus S_{i}$, is
of standard structure, if $u$ is adjacent to every vertex in
$S_{i}\cap Y$ and $v$ is adjacent to every vertex in $S_{i}\cap X$,
particularly, we denote the tree by $T_{u, v}$. Denote the set of
trees in $\mathfrak{A}_{0}$ with the standard structure by
$\mathcal{A}_{0}$, clearly, $\mathcal{A}_{0}=\mathfrak{A}_{0}$.
Similarly, denote the set of trees in $\mathfrak{A}_{1}$ and
$\mathfrak{A}_{2}$ with the standard structure by $\mathcal{A}_{0}$
and $\mathcal{A}_{2}$, respectively.

\begin{lem}\label{lem2}
Let $A$ be a maximum set of internally disjoint trees connecting
$S_{i}$, $A \subset
\mathfrak{A}_{0}\cup\mathfrak{A}_{1}\cup\mathfrak{A}_{2}$. Then we
can always find a set $A''$ of internally disjoint trees connecting
$S_{i}$, such that $\mid A \mid = \mid A'' \mid$ and $A'' \subset
\mathcal{A}_{0}\cup\mathcal{A}_{1}\cup\mathcal{A}_{2}$.
\end{lem}
\pf Suppose there is a tree $T^{0}$ in $A$ such that $T^{0}\in
\mathfrak{A}_{1}$ but $T^{0} \notin \mathcal{A}_{1}$, and $V(T^{0})=
S_{i}\cup \{u_{0}\}$, where $u_{0}\in X\setminus S_{i}$. Note that
the case $u_{0}\in Y\setminus S_{i}$ is similar. Since $T^{0} \notin
\mathcal{A}_{1}$, there are some vertices in $S_{i}\cap Y$, say
$y_{i_{1}}, \ldots , y_{i_{t}}$, not adjacent to $u_{0}$. Then we
can connect $y_{i_{1}}$ to $u_{0}$ by a new edge. It will produce a
unique cycle. Delete the other edge incident with $y_{i_{1}}$ on the
cycle. The graph remains a tree. Do the operation to $y_{i_{2}},
\ldots , y_{i_{t}}$ in turn. Finally we get a tree $T$ whose vertex
set is $S_{i} \cup \{u_{0}\}$ and $u_{0}$ is adjacent to every
vertex in $S_{i}\cap Y$, that is, $T$ is of standard structure. For
each tree $T_{n}\in A\setminus\{T^{0}\}$, clearly $T_{n}$ does not
contain $u_{0}$, including the edges incident with $u_{0}$. So
$V(T)\cap V(T_{n}) = S_{i}$ and $E(T)\cap E(T_{n}) =\emptyset$.
Replace $A$ by $A^{1}=A\setminus\{T^{0}\}\cup\{T\}$, and then
$A^{1}$ contains less trees not in $\mathcal{A}_{0}\cup
\mathcal{A}_{1}\cup \mathcal{A}_{2}$ than $A$. Suppose that there is
a tree $T^{1}$ in $A$ such that $T^{1}\in \mathfrak{A}_{2}$ but
$T^{1} \notin \mathcal{A}_{2}$ and $V(T^{1})= S_{i}\cup \{u_{1},
v_{1}\}$, where $u_{1}\in X\setminus S_{i}$ and $v_{1}\in Y\setminus
S_{i}$. $T_{u_{1}, v_{1}}$ is the tree in $\mathcal{A}_{2}$ whose
vertex set is $S_{i}\cup \{u_{1}, v_{1}\}$. Then for each tree
$T_{n}\in A\setminus\{T^{1}\}$, $V(T_{u_{1}, v_{1}})\cap V(T_{n}) =
S_{i}$ and $E(T_{u_{1}, v_{1}})\cap E(T_{n}) = \emptyset$. Replace
$A$ by $A^{1}=A\setminus\{T^{1}\}\cup\{T_{u_{1}, v_{1}}\}$. Then
$A^{1}$ contains less trees not in $\mathcal{A}_{0}\cup
\mathcal{A}_{1}\cup \mathcal{A}_{2}$ than $A$. Repeating the
process, we can get a series of sets $A^{0}, A^{1}, \ldots, A^{t}$,
such that $A^{0}=A$, $A^{t}=A''$, and $A^{j}$ contains less trees
not in $\mathcal{A}_{0}\cup \mathcal{A}_{1}\cup \mathcal{A}_{2}$
than $A^{j-1}$, for $1\leq j \leq t$, where all $A^{s}$ are sets of
internally disjoint trees connecting $S_{i}$, $A^{s} \subset
\mathfrak{A}_{0}\cup\mathfrak{A}_{1}\cup\mathfrak{A}_{2}$, for
$0\leq s \leq t$, and $\mid A^{0}\mid=\cdots=\mid A^{t}\mid$. So we
finally get the set $A'' \subset
\mathcal{A}_{0}\cup\mathcal{A}_{1}\cup\mathcal{A}_{2}$ which has the
same cardinality as $A$. \qed

So, we can assume that the maximum set $A$ of internally disjoint
trees connecting $S_{i}$ is contained in
$\mathcal{A}_{0}\cup\mathcal{A}_{1}\cup\mathcal{A}_{2}$. Namely, all
trees in $A$ are of standard structure.

For simplicity, we denote the union of the vertex sets of all trees
in set $A$ by $V(A)$ and the union of the edge sets of all trees in
set $A$ by $E(A)$. Let $A$ be a set of internally disjoint trees
connecting $S_{i}$. Let $A_{0}:= A\cap \mathcal{A}_{0}$, $A_{1}:=
A\cap \mathcal{A}_{1}$ and $A_{2}:= A\cap \mathcal{A}_{2}$. Then
$A=A_{0}\cup A_{1}\cup A_{2}$. Let $U(A):=V(G)\setminus V(A)$.

\begin{lem}\label{lem3}
Let $A\subset \mathcal{A}_{0}\cup\mathcal{A}_{1}\cup\mathcal{A}_{2}$
be a maximum set of internally disjoint trees connecting $S_{i}$,
$A=A_{0}\cup A_{1}\cup A_{2}$ and $U(A):=V(G)\setminus V(A)$. Then
either $U(A)\cap X=\emptyset$ or $U(A)\cap Y=\emptyset$.
\end{lem}
\pf If $U(A)\cap X\neq\emptyset$ and $U(A)\cap Y\neq\emptyset$, let
$x\in U(A)\cap X$ and $y\in U(A)\cap Y$. Then the tree $T_{x,y} \in
\mathcal{A}_{2}$ with vertex set $S_{i}\cup \{x, y\}$ is a tree that
connects $S_{i}$. Moreover, $V(T)\cap V(A)=S_{i}$ and for any tree
$T' \in A$, $T$ and $T'$ are edge-disjoint. So, $A\cup \{T\}$ is
also a set of internally disjoint trees connecting $S_{i}$,
contradicting to the maximality of $A$. \qed

So we conclude that if $A$ is a maximum set of internally disjoint
trees connecting $S_{i}$, then $U(A)\subset X$ or $U(A)\subset Y$.

\begin{lem}\label{lem4}
Let $A\subset \mathcal{A}_{0}\cup\mathcal{A}_{1}\cup\mathcal{A}_{2}$
be a maximum set of internally disjoint trees connecting $S_{i}$,
$A=A_{0}\cup A_{1}\cup A_{2}$ and $U(A):=V(G)\setminus V(A)$. If
$U(A)\neq \emptyset$ and $A_{0}\neq \emptyset$, then we can find a
set $A'=A'_{0}\cup A'_{1}\cup A'_{2}$ of internally disjoint trees
connecting $S_{i}$, such that $|A'_{0}|=|A_{0}|-1$,
$|A'_{1}|=|A_{1}|+1$, $A'_{2}=A_{2}$ and $|U(A')|=|U(A)|-1$.
\end{lem}
\pf Let $u \in U(A)$ and $T \in A_{0}$. Without loss of generality,
suppose $u\in X$. Then we can connect $u$ to $y_{1}$ by a new edge,
and the new graph becomes a tree $T' \in \mathfrak{A}_{1}$. Using
the method in Lemma \ref{lem2}, we can transform $T'$ into a tree
$T''$ with the standard structure. Then $T'' \in \mathcal{A}_{1}$.
Let $A'_{0}:=A_{0}\setminus T$, $A'_{1}:=A_{1} \cup \{T''\}$ and
$A'_{2}=A_{2}$. It is easy to see that $A'=A'_{0}\cup A'_{1}\cup
A'_{2}$ is a set of internally disjoint trees connecting $S_{i}$.
Since $|A'_{0}|=|A_{0}|-1$, $|A'_{1}|=|A_{1}|+1$, and
$A'_{2}=A_{2}$, $A'$ is a maximum set of internally disjoint trees
connecting $S_{i}$ and $|U(A')|=|U(A)|-1$. \qed

So, we can assume that for the maximum set $A$ of internally
disjoint trees connecting $S_{i}$, either $U(A)=\emptyset$ or
$A_{0}= \emptyset$. Moreover, if $A'$ is a set of internally
disjoint trees connecting $S_{i}$ which we find currently,
$U(A')\neq \emptyset$ and the edges in $E(G[S_{i}])\setminus E(A')$
can form a tree $T$ in $\mathcal{A}_{0}$, then we will add to $A'$
the tree $T''$ in Lemma \ref{lem4} rather than the tree $T$.

\begin{lem}\label{lem5}
Let $A\subset \mathcal{A}_{0}\cup\mathcal{A}_{1}\cup\mathcal{A}_{2}$
be a maximum set of internally disjoint trees connecting $S_{i}$,
$A=A_{0}\cup A_{1}\cup A_{2}$ and $U(A):=V(G)\setminus V(A)$. If
there is a vertex $x \in U(A) \subset X$ and a tree $T \in A_{1}$
with vertex set $S_{i} \cup \{y\}$, where $y \in Y\setminus S_{i}$.
Then we can find a set $A'=A'_{0}\cup A'_{1}\cup A'_{2}$ of
internally disjoint trees connecting $S_{i}$, such that
$A'_{0}=A_{0}$, $|A'_{1}|=|A_{1}|-1$, $|A'_{2}|=|A_{2}|+1$ and
$|U(A')|=|U(A)|-1$.
\end{lem}
\pf Let $T_{x,y}$ be the tree in $\mathcal{A}_{2}$ whose vertex set
is $S_{i}\cup \{x,y\}$. Then $A'=A\setminus T \cup \{T_{x,y}\}$ is
just the set we want. \qed

The case that there is a vertex $y \in U(A) \subset Y$ and a tree $T
\in A_{1}$ with vertex set $S_{i} \cup \{x\}$, where $x \in
X\setminus S_{i}$, is similar. So we can assume that, for the
maximum set $A$ of internally disjoint trees connecting $S_{i}$, $A$
satisfies one of the following properties:

(1) $U(A)=\emptyset$

(2) $\emptyset\neq U(A)\subset X$ and $V(A_{1})\setminus S_{i}\subset X$

(3) $\emptyset\neq U(A)\subset Y$ and $V(A_{1})\setminus S_{i}\subset Y$

Now, we can see that if $U(A)\neq \emptyset$, then all vertices in
$V(A_{1})\setminus S_{i}$ belong to the same partition. Next, we
will show that we can always find a set $A$ of internally disjoint
trees connecting $S_{i}$, such that no matter whether $U(A)$ is
empty, all vertices in $V(A_{1})\setminus S_{i}$ belong to the same
partition. To show this, we need the following lemma.

\begin{lem}\label{lem6}
Let $p$, $q$ be two nonnegative integers. If $p(k-1)+qi\leq i(k-i)$,
and there are $q$ vertices $u_{1}, u_{2}, \ldots, u_{q}\in
X\setminus S_{i}$, then we can always find $p$ trees $T_{1}, T_{2},
\ldots, T_{p}$ in $\mathcal{A}_{0}$ and $q$ trees $T_{p+1}, T_{p+2},
\ldots, T_{p+q}$ in $\mathcal{A}_{1}$, such that $V(T_{j})=S_{i}$
for $1\leq j \leq p$, $V(T_{p+m})=S_{i}\cup\{u_{m}\}$ for $1\leq m
\leq q$, and $T_{r}$ and $T_{s}$ are edge-disjoint for $1\leq r < s
\leq p+q$. Similarly, if $p(k-1)+q(k-i)\leq i(k-i)$, and there are
$q$ vertices $v_{1}, v_{2}, \ldots, v_{q}\in Y\setminus S_{i}$, then
we can always find $p$ trees $T_{1}, T_{2}, \ldots, T_{p}$ in
$\mathcal{A}_{0}$ and $q$ trees $T_{p+1}, T_{p+2}, \ldots, T_{p+q}$
in $\mathcal{A}_{1}$, such that $V(T_{j})=S_{i}$ for $1\leq j \leq
p$, $V(T_{p+m})=S_{i}\cup\{v_{m}\}$ for $1\leq m \leq q$, and
$T_{r}$ and $T_{s}$ are edge-disjoint for $1\leq r < s \leq p+q$.
\end{lem}
\pf If $p(k-1)+qi\leq i(k-i)$, then $p(k-1)\leq i(k-i)$, namely $p
\leq \lfloor\frac{i(k-i)}{k-1}\rfloor$. Then with the method which
we used to find edge-disjoint spanning trees in the proof of Theorem
\ref{thm1}, we can find $p$ edge-disjoint trees $T_{1}, T_{2},
\ldots, T_{p}$ in $\mathcal{A}_{0}$, just by taking $a=i$, $b=k-i$
and $t=p$. Moreover, let $D_{s}^{p}$ denote the number of edges
incident with $x_{s}$ in all of the $p$ trees, then according to the
method, $|D_{s}^{p}-D_{t}^{p}|\leq 1$ for $1\leq s, t \leq i$. Now,
denote by $B_{s}^{p}$ the number of edges incident with $x_{s}$
which we have not used in the $p$ trees. Then
$|B_{s}^{p}-B_{t}^{p}|\leq 1$ for $1\leq s, t \leq i$. Since
$B_{1}^{p}+B_{2}^{p}+\cdots +B_{i}^{p}=i(k-i)-p(k-1)\geq qi$,
$B_{s}^{p} \geq q$. Because for each tree in $\mathcal{A}_{1}$ with
vertex set $S_{i}\cup\{u\}$, where $u\in X\setminus S_{i}$, the
vertices in $S_{i}\cap X$ all have degree $1$, we can find $q$
edge-disjoint trees $T_{p+1}, T_{p+2}, \ldots, T_{p+q}$ in
$\mathcal{A}_{1}$. Since the edges in $T_{p+1}, T_{p+2}, \ldots,
T_{p+q}$ are not used in $T_{1}, T_{2}, \ldots, T_{p}$ for $1\leq r
< s \leq p+q$, $T_{r}$ and $T_{s}$ are edge-disjoint. The proof of
the second half of the lemma is similar. \qed

\begin{lem}\label{lem7}
Let $A\subset \mathcal{A}_{0}\cup\mathcal{A}_{1}\cup\mathcal{A}_{2}$
be a maximum set of internally disjoint trees connecting $S_{i}$,
$A=A_{0}\cup A_{1}\cup A_{2}$ and $U(A):=V(G)\setminus V(A)$. If
there are $s$ trees $T_{1}, T_{2}, \ldots, T_{s} \in A_{1}$ with
vertex set $S_{i} \cup \{u^{1}\}$, $S_{i} \cup \{u^{2}\}$, $\ldots$,
$S_{i} \cup \{u^{s}\}$ respectively, where $u^{j} \in X\setminus
S_{i}$ for $1\leq j \leq s$, and $t$ trees $T_{s+1}, T_{s+2},
\ldots, T_{s+t} \in A_{1}$ with vertex set $S_{i} \cup \{v^{1}\}$,
$S_{i} \cup \{v^{2}\}$, $\ldots$, $S_{i} \cup \{v^{t}\}$
respectively, where $v^{j} \in Y\setminus S_{i}$ for $1\leq j \leq
t$. Then we can find a set $A'=A'_{0}\cup A'_{1}\cup A'_{2}$ of
internally disjoint trees connecting $S_{i}$, such that $|A|=|A'|$
and all vertices in $V(A'_{1})\setminus S_{i}$ belong to the same
partition.
\end{lem}
\pf Let $|A_{0}|=p$. Since $A$ is a set of internally disjoint trees
connecting $S_{i}$, we have $p(k-1)+si+t(k-i)\leq i(k-i)$, where
$si$ denote the $si$ edges incident with $x_{1}, \ldots, x_{i}$ in
$T_{1}, T_{2}, \ldots, T_{s}$, and $t(k-i)$ denote the $t(k-i)$
edges incident with $y_{1}, \ldots, y_{k-i}$ in $T_{s+1}, T_{s+2},
\ldots, T_{s+t}$. If $s\leq t$, then
$p(k-1)+si+s(k-i)+(t-s)(k-i)\leq i(k-i)$, and hence
$(p+s)(k-1)+(t-s)(k-i)\leq i(k-i)$. Obviously, there are $t-s$
vertices $v^{s+1}, v^{s+2}, \ldots, v^{t}\in Y\setminus S_{i}$, and
therefore by Lemma \ref{lem6}, we can find $p+s$ trees in
$\mathcal{A}_{0}$ and $t-s$ trees in $\mathcal{A}_{1}$, such that
all these trees are internally disjoint trees connecting $S_{i}$.
Now let $A'_{0}$ be the set of the $p+s$ trees in $\mathcal{A}_{0}$,
$A'_{1}$ be the set of the $t-s$ trees in $\mathcal{A}_{1}$ and
$A'_{2}:=A_{2}\cup \{T_{u^{j}, v^{j}}, 1\leq j \leq s\}$. Then
$A'=A'_{0}\cup A'_{1}\cup A'_{2}$ is just the set we want. The case
that $s > t$ is similar. \qed

From Lemmas \ref{lem5} and \ref{lem7}, we can see that, if $A'$ is a
set of internally disjoint trees connecting $S_{i}$ which we find
currently, $U(A')\cap X\neq \emptyset$ and $U(A')\cap Y\neq
\emptyset$, then no matter how many edges are there in
$E(G[S_{i}])\setminus E(A')$, we always add to $A'$ the trees in
$\mathcal{A}_{2}$ rather than the trees in $\mathcal{A}_{1}$.

Next, let us state and prove our main result.

\begin{thm}\label{thm2}
Given any two positive integers $a$ and $b$, let $K_{a, b}$ denote a
complete bipartite graph with a bipartition of sizes $a$ and $b$,
respectively. Then we have the following results: if $k>b-a+2$ and
$a-b+k$ is odd then
$$
\kappa_{k}(K_{a,b})=\frac{a+b-k+1}{2}+\lfloor\frac{(a-b+k-1)(b-a+k-1)}{4(k-1)}\rfloor;
$$
if $k>b-a+2$ and $a-b+k$ is even then
$$
\kappa_{k}(K_{a,b})=\frac{a+b-k}{2}+\lfloor\frac{(a-b+k)(b-a+k)}{4(k-1)}\rfloor;
$$
and if $k\leq b-a+2$ then
$$
\kappa_{k}(K_{a,b})=a.
$$
\end{thm}
\pf Recall that $\kappa_{k}(G) = \min\{\kappa(S)\}$, where the
minimum is taken over all $k$-element subsets $S$ of $V(G)$. Let
$X=\{x_{1}, x_{2}, \ldots, x_{a}\}$ and $Y=\{y_{1}, y_{2}, \ldots,
y_{b}\}$ be the bipartition of $K_{a,b}$, where $a\leq b$. As we
have mentioned, all vertices in $X$ are equivalent and all vertices
in $Y$ are equivalent. So instead of considering all $k$-element
subsets $S$ of $V(G)$, we can restrict our attention to the subsets
$S_{i}$, for $0 \leq i \leq k$, where $S_{i}$ is an $k$-element
subsets of $V(G)$ such that $S_{i}\cap X = \{x_{1}, x_{2}, \ldots,
x_{i}\}$, $S_{i}\cap Y = \{y_{1}, y_{2}, \ldots, y_{k-i}\}$, $1 \leq
i \leq k$ and $S_{0}\cap X =\emptyset$, $S_{0}\cap Y =\{y_{1},
y_{2}, \ldots, y_{k}\}$. Notice that, if $i>a$ or $k-i > b$ then
$S_{i}$ does not exist, and if $k > b$ then $S_{0}$ does not exist.
So, we need only to consider $S_{i}$ for $\max\{0,\ k-b\} \leq i
\leq \min\{a,\ k\}$.

From the above lemmas, we can decide our principle to find the
maximum set of internally disjoint trees connecting $S_{i}$. Namely,
first we find as many trees in $\mathcal{A}_{2}$ as possible, next
we find as many trees in $\mathcal{A}_{1}$ as possible, and finally
we find as many trees in $\mathcal{A}_{0}$ as possible.

For a set $S_{i}=\{x_{1}, x_{2}, \ldots, x_{i}, y_{1}, y_{2},
\ldots, y_{k-i}\}$, let $A$ be the maximum set of internally
disjoint trees connecting $S_{i}$ we find with our principle. We now
compute $|A|$.

\textbf{Case 1.} $k\leq b-a+2$

Obviously, $\kappa(S_{0})=a$.

For $S_{1}$, since $k\leq b-a+2$, then
$$
b-(k-1)=b-k+1\geq a-2+1=a-1.
$$
So, $|A_{2}|=a-1$. If $b-k+1=a-1$, then $|A_{1}|=0$, $|A_{0}|=1$. If
$b-k+1>a-1$, then $|A_{1}|=1$, $|A_{0}|=0$. No matter which case
happens, we have $\kappa(S_{1})=|A_{2}|+|A_{1}|+|A_{0}|=a$.

For $S_{i}$, $i\geq 2$, since $k\leq b-a+2$, then
$$
b-(k-i)=b-k+i\geq a-2+i>a-i.
$$
So, $|A_{2}|=a-i$. Since $b-k+i-(a-i)=b-a-k+2i\geq -2+2i\geq i$,
then $|A_{1}|=i$ and $|A_{0}|=0$. Thus
$\kappa(S_{i})=|A_{2}|+|A_{1}|+|A_{0}|=a$.

In summary, if $k\leq b-a+2$, then $\kappa_{k}(G)=a$.

\textbf{Case 2.} $k > b-a+2$

First, let us compare $\kappa(S_{i})$ with $\kappa(S_{k-i})$,
for $0\leq i \leq \lfloor\frac{k}{2}\rfloor$. If $a=b$,
clearly, $\kappa(S_{i})=\kappa(S_{k-i})$. So we may assume that $a<b$.

For $i=0$, $\kappa(S_{0})=a < b = \kappa(S_{k})$.

For $1 \leq i \leq \lfloor\frac{k}{2}\rfloor$, we will give the
expressions of $\kappa(S_{i})$ and $\kappa(S_{k-i})$.

First for $S_i$, since every pair of vertices $u \in X\setminus S_i$
and $v \in Y\setminus S_i$ can form a tree $T_{u,v}$, then
$|A_2|=\min \{a-i, b-(k-i)\}$. Namely,
\[
|A_2|=\left\{
\begin{array}{ll}
a-i &\mbox{if ~$i \geq \frac{a-b+k}{2}$~;}\\
b-k+i &\mbox{if ~$i < \frac{a-b+k}{2}$~.}
\end{array}
\right.
\]
Next, since every tree $T$ in $A_1$ has a vertex in $V\setminus (S_i\cup V(A_2))$,
we have
\[
|A_1|\leq\left\{
\begin{array}{ll}
b-k+i-(a-i) &\mbox{if ~$i \geq \frac{a-b+k}{2}$~;}\\
a-i-(b-k+i) &\mbox{if ~$i < \frac{a-b+k}{2}$~.}
\end{array}
\right.
\]
On the other hand, if the tree $T$ has vertex set $S_i\cup\{u\}$,
where $u \in X\setminus S_i$, then every vertex in $S_i\cap X$ is
incident with one edge in $E(S_i)$, where $E(S_i)$ denotes the set
of edges whose ends are both in $S_i$. And if the tree $T$ has
vertex set $S_i\cup\{v\}$, where $v \in Y\setminus S_i$, then every
vertex in $S_i\cap Y$ is incident with one edge in $E(S_i)$. Since
every vertex in $S_i\cap X$ is incident with $k-i$ edges in $E(S_i)$
and every vertex in $S_i\cap Y$ is incident with $i$ edges in
$E(S_i)$, we have
\[
|A_1|\leq\left\{
\begin{array}{ll}
i &\mbox{if ~$i \geq \frac{a-b+k}{2}$~;}\\
k-i &\mbox{if ~$i < \frac{a-b+k}{2}$~.}
\end{array}
\right.
\]
Combining the two inequalities, we get
\[
|A_1|=\left\{
\begin{array}{ll}
\min\{b-a-k+2i, i\} &\mbox{if ~$i \geq \frac{a-b+k}{2}$~;}\\
\min\{a-b+k-2i, k-i\} &\mbox{if ~$i < \frac{a-b+k}{2}$~.}
\end{array}
\right.
\]
Thus
\[
|A_1|=\left\{
\begin{array}{ll}
i &\mbox{if ~$i \geq a-b+k$~;}\\
b-a-k+2i &\mbox{if ~$\frac{a-b+k}{2}\leq i < a-b+k $~;}\\
a-b+k-2i &\mbox{if ~$i < \frac{a-b+k}{2}$~.}
\end{array}
\right.
\]
Finally, by Lemma \ref{lem6} we have
\[
|A_0|=\left\{
\begin{array}{ll}
\lfloor\frac{i(k-i)-|A_1|(k-i)}{k-1}\rfloor &\mbox{if ~$i \geq \frac{a-b+k}{2}$~;}\\
\lfloor\frac{i(k-i)-|A_1|i}{k-1}\rfloor &\mbox{if ~$i < \frac{a-b+k}{2}$~.}
\end{array}
\right.
\]
Thus
\[
|A_0|=\left\{
\begin{array}{ll}
0 &\mbox{if ~$i \geq a-b+k$~;}\\
\lfloor\frac{[i-(b-a-k+2i)](k-i)}{k-1}\rfloor &\mbox{if ~$\frac{a-b+k}{2}\leq i < a-b+k $~;}\\
\lfloor\frac{[k-i-(a-b+k-2i)]i}{k-1}\rfloor &\mbox{if ~$i < \frac{a-b+k}{2}$~.}
\end{array}
\right.
\]
And hence
\[
\kappa(S_i)=\left\{
\begin{array}{ll}
a &\mbox{if ~$i \geq a-b+k$~;}\\
b-k+i+\lfloor\frac{[i-(b-a-k+2i)](k-i)}{k-1}\rfloor &\mbox{if ~$\frac{a-b+k}{2}\leq i < a-b+k $~;}\\
a-i+\lfloor\frac{[k-i-(a-b+k-2i)]i}{k-1}\rfloor &\mbox{if ~$i < \frac{a-b+k}{2}$~.}
\end{array}
\right.
\]
Notice that $i\geq 1$, and hence $k-i \leq k-1$.

If $\frac{a-b+k}{2}\leq i < a-b+k $, then
$$
\lfloor\frac{[i-(b-a-k+2i)](k-i)}{k-1}\rfloor\leq i-(b-a-k+2i)= a-b+k-i.
$$
So, $\kappa(S_i)\leq b-k+i+a-b+k-i=a$.

If $i < \frac{a-b+k}{2}$, then $a-b+k-2i>0$, $k-i-(a-b+k-2i)<k-i\leq k-1$, and hence
$$
\lfloor\frac{[k-i-(a-b+k-2i)]i}{k-1}\rfloor\leq i.
$$
So, $\kappa(S_i)\leq a-i+i=a$

Thus $\kappa(S_i)\leq a$, for $i\geq 1$.

Next, considering $S_{k-i}$, similarly, we have
$$
|A_2|=\min\{a-(k-i), b-i\}.
$$
Since $a<b$ and $i\leq \lfloor\frac{k}{2}\rfloor \leq
\lceil\frac{k}{2}\rceil \leq k-i$, then $b-i > a-(k-i)$. So $|A_2|=
a-k+i$ and $|A_1|=\min \{b-i-(a-k+i), k-i\}$. Hence
\[
|A_1|=\left\{
\begin{array}{ll}
k-i &\mbox{if ~$i \leq b-a$~;}\\
b-a+k-2i &\mbox{if ~$i > b-a$~.}
\end{array}
\right.
\]
Moreover,
\[
|A_0|=\left\{
\begin{array}{ll}
0 &\mbox{if ~$i \leq b-a$~;}\\
\lfloor\frac{[k-i-(b-a+k-2i)]i}{k-1}\rfloor &\mbox{if ~$i > b-a$~.}
\end{array}
\right.
\]
So,
\[
\kappa(S_{k-i})=\left\{
\begin{array}{ll}
a &\mbox{if ~$i \leq b-a$~;}\\
b-i+\lfloor\frac{[k-i-(b-a+k-2i)]i}{k-1}\rfloor &\mbox{if ~$i > b-a$~.}
\end{array}
\right.
\]
Now, we can compare $\kappa(S_{i})$ with $\kappa(S_{k-i})$.
For $i\leq b-a$, $\kappa(S_{k-i})=a\geq \kappa(S_i)$. For
$i > b-a$, there must be $b-a< k-i$, that is, $i< a-b+k$.

If $\frac{a-b+k}{2}\leq i < a-b+k $, then
\begin{eqnarray*}
\kappa(S_{k-i})-\kappa(S_i)&=& b-i+\lfloor\frac{[k-i-(b-a+k-2i)]i}{k-1}\rfloor\\
&&-\{b-k+i+\lfloor\frac{[i-(b-a-k+2i)](k-i)}{k-1}\rfloor\}\\
&\geq &(k-2i)+\lfloor\frac{(k-2i)(b-a-k)}{k-1}\rfloor\\
&\geq&(k-2i)+\lfloor\frac{(k-2i)(1-k)}{k-1}\rfloor\\
&\geq&(k-2i)-(k-2i)=0.
\end{eqnarray*}
So, $\kappa(S_{k-i})\geq\kappa(S_i)$.

If $i < \frac{a-b+k}{2}$, then
\begin{eqnarray*}
\kappa(S_{k-i})-\kappa(S_i)&=& b-i+\lfloor\frac{[k-i-(b-a+k-2i)]i}{k-1}\rfloor\\
&&-\{a-i+\lfloor\frac{[k-i-(a-b+k-2i)]i}{k-1}\rfloor\}\\
&\geq&(b-a)+\lfloor\frac{(2i)(a-b)}{k-1}\rfloor.
\end{eqnarray*}
Since $i < \frac{a-b+k}{2}$, then $2i\leq k-1$, and hence
$\frac{(2i)(a-b)}{k-1}\geq a-b$. So,
$\kappa(S_{k-i})-\kappa(S_i)\geq b-a+a-b=0$. Thus,
$\kappa(S_{k-i})\geq\kappa(S_i)$.

In summary, $\kappa(S_{k-i})\geq\kappa(S_i)$, for $0\leq i \leq
\lfloor\frac{k}{2}\rfloor$. So, in order to get $\kappa_{k}(G)$, it
is enough to consider $\kappa(S_{i})$, for $0\leq i \leq
\lfloor\frac{k}{2}\rfloor$.

Next, let us compare $\kappa(S_{i})$ with $\kappa(S_{i+1})$,
for $0\leq i \leq \lfloor\frac{k}{2}\rfloor-1$. For
$i=0$, $\kappa(S_{i})=a\geq \kappa(S_{i+1})$. For $1\leq i \leq \lfloor\frac{k}{2}\rfloor-1$,
\[
\kappa(S_i)=\left\{
\begin{array}{ll}
a &\mbox{if ~$i \geq a-b+k$~;}\\
b-k+i+\lfloor\frac{[i-(b-a-k+2i)](k-i)}{k-1}\rfloor &\mbox{if ~$\frac{a-b+k}{2}\leq i < a-b+k $~;}\\
a-i+\lfloor\frac{[k-i-(a-b+k-2i)]i}{k-1}\rfloor &\mbox{if ~$i < \frac{a-b+k}{2}$~.}
\end{array}
\right.
\]
and
\[
\kappa(S_{i+1})=\left\{
\begin{array}{ll}
a &\mbox{if ~$i \geq a-b+k-1$~;}\\
b-k+i+1+\lfloor\frac{[i+1-(b-a-k+2i+2)](k-i-1)}{k-1}\rfloor &\mbox{if ~$\frac{a-b+k}{2}-1\leq i < a-b+k-1 $~;}\\
a-i-1+\lfloor\frac{[k-i-1-(a-b+k-2i-2)](i+1)}{k-1}\rfloor &\mbox{if ~$i < \frac{a-b+k}{2}-1$~.}
\end{array}
\right.
\]
So, $\kappa(S_{a-b+k})=\kappa(S_{a-b+k+1})=\cdots=\kappa(S_{\min\{a,k\}})=a$.

If $i<\frac{a-b+k}{2}-1$, then
\begin{eqnarray*}
\kappa(S_{i})-\kappa(S_{i+1})&=& a-i+\lfloor\frac{[k-i-(a-b+k-2i)]i}{k-1}\rfloor\\
&&-\{a-i-1+\lfloor\frac{[k-i-1-(a-b+k-2i-2)]i+1}{k-1}\rfloor\}\\
&\geq&1+\lfloor\frac{(a-b-2i-1)}{k-1}\rfloor\\
&\geq&1+\lfloor\frac{1-k}{k-1}\rfloor\\
&\geq&1-1=0.
\end{eqnarray*}
So, $\kappa(S_{i})\geq\kappa(S_{i+1})$. Namely, if $a-b+k$ is odd,
we have
$$
\kappa(S_{0})\geq\kappa(S_{1})\geq\cdots\geq\kappa(S_{\frac{a-b+k-3}{2}})\geq\kappa(S_{\frac{a-b+k-1}{2}}).
$$
and if $a-b+k$ is even, we have
$$
\kappa(S_{0})\geq\kappa(S_{1})\geq\cdots\geq\kappa(S_{\frac{a-b+k-4}{2}})\geq\kappa(S_{\frac{a-b+k-2}{2}}).
$$

If $i=\frac{a-b+k}{2}-1$, $\kappa(S_i)=\frac{a+b-k}{2}+1+\lfloor\frac{(b-a+k-2)(a-b+k-2)}{4(k-1)}\rfloor$.

If $i=\frac{a-b+k-1}{2}$, $\kappa(S_i)=\frac{a+b-k+1}{2}+\lfloor\frac{(b-a+k-1)(a-b+k-1)}{4(k-1)}\rfloor$.

If $i=\frac{a-b+k}{2}$, $\kappa(S_i)=\frac{a+b-k}{2}+\lfloor\frac{(b-a+k)(a-b+k)}{4(k-1)}\rfloor$.

If $i=\frac{a-b+k+1}{2}$, $\kappa(S_i)=\frac{a+b-k+1}{2}+\lfloor\frac{(b-a+k-1)(a-b+k-1)}{4(k-1)}\rfloor$.

If $a-b+k$ is even, since
\begin{eqnarray*}
&&(a-b+k)(b-a+k)-(b-a+k-2)(a-b+k-2)\\
&=&(a-b+k)(b-a+k)-[(a-b+k)(b-a+k)-2(b-a+k)-2(a-b+k-2)]\\
&=&4(k-1),
\end{eqnarray*}
then we have $\kappa(S_{\frac{a-b+k}{2}-1})=\kappa(S_{\frac{a-b+k}{2}})$. If $a-b+k$ is odd, we have
$\kappa(S_{\frac{a-b+k-1}{2}})=\kappa(S_{\frac{a-b+k+1}{2}})$.

If $\frac{a-b+k}{2}\leq i \leq a-b+k-1 $, then
\begin{eqnarray*}
\kappa(S_{i+1})-\kappa(S_{i})&=& b-k+i+1+\lfloor\frac{[i+1-(b-a-k+2i+2)](k-i-1)}{k-1}\rfloor\\
&&-\{b-k+i+\lfloor\frac{[i-(b-a-k+2i)](k-i)}{k-1}\rfloor\}\\
&\geq&1+\lfloor\frac{(b-a-2k+2i+1)}{k-1}\rfloor\\
&\geq&1+\lfloor\frac{1-k}{k-1}\rfloor\\
&\geq&1-1=0.
\end{eqnarray*}
So, $\kappa(S_{i+1})\geq\kappa(S_i)$. Namely, if $a-b+k$ is odd, we
have
$$
\kappa(S_{\frac{a-b+k+1}{2}})\leq\kappa(S_{\frac{a-b+k+3}{2}})\leq\cdots\leq\kappa(S_{a-b+k-1})\leq\kappa(S_{a-b+k})=a,
$$
and if $a-b+k$ is even, we have
$$
\kappa(S_{\frac{a-b+k}{2}})\leq\kappa(S_{\frac{a-b+k+2}{2}})\leq\cdots\leq\kappa(S_{a-b+k-1})\leq\kappa(S_{a-b+k})=a.
$$
Thus, if $k>b-a+2$ and $a-b+k$ is odd,
$$
\kappa_{k}(K_{a,b})=\kappa(S_{\frac{a-b+k-1}{2}})=\frac{a+b-k+1}{2}+\lfloor\frac{(a-b+k-1)(b-a+k-1)}{4(k-1)}\rfloor,
$$
and if $k>b-a+2$ and $a-b+k$ is even,
$$
\kappa_{k}(K_{a,b})=\kappa(S_{\frac{a-b+k}{2}})=\frac{a+b-k}{2}+\lfloor\frac{(a-b+k)(b-a+k)}{4(k-1)}\rfloor.
$$
The proof is complete.\qed

Notice that, when $k=a+b$, the result coincides with Theorem
\ref{thm1}.

\end{document}